 \documentclass{amsart}

\usepackage{amsmath,amsthm,amsfonts, amssymb,amscd}
\usepackage[all]{xy}

\newtheorem{theorem}{Theorem}[section]
\newtheorem{lemma}[theorem]{Lemma}%[section]
\newtheorem{example}[theorem]{Example}%[section]
\newtheorem{remark}[theorem]{Remark}%[section]
\newtheorem{corollary}[theorem]{Corollary}%[section]
\newtheorem{proposition}[theorem]{Proposition}

\newcommand{\dx}{\, \mbox{\rm d}}

\newcommand{\Aff}{\mbox{\rm Aff}}

\begin{document}
\title[Representation of States on Effect-Tribes and Effect Algebras]{Representation
of States on Effect-Tribes and Effect Algebras by Integrals}
\author{Anatolij Dvure\v censkij}
\date{}%Jan. 5, 2010
\maketitle

\begin{center}
\footnote{


The  author thanks  for the support by Center of Excellence SAS
-~Quantum Technologies~-,  ERDF OP R\&D Projects CE QUTE ITMS
26240120009 and meta-QUTE ITMS 26240120022, the grant VEGA No.
2/0032/09 SAV. }
\small{Mathematical Institute,  Slovak Academy of Sciences\\
\v Stef\'anikova 49, SK-814 73 Bratislava, Slovakia\\
E-mail: {\tt dvurecen@mat.savba.sk}, }
\end{center}

\begin{abstract}  We  describe $\sigma$-additive states on
effect-tribes  by integrals. Effect-tribes are monotone
$\sigma$-complete effect algebras of functions where operations are
defined by points. Then we show that every state on an effect
algebra is an integral through a Borel regular probability measure.
Finally, we show that every $\sigma$-convex combination of extremal
states on a monotone $\sigma$-complete effect algebra is a
Jauch-Piron state.
\end{abstract}

\section{Introduction} %1

The study of the mathematical foundations of quantum mechanics was
initiated by the famous article by Birkhoff and von Neumann
\cite{BiNe} and the theory was called quantum logics or in the
present the theory of quantum structures. In 1993 the theory of
quantum structures was enriched by {\it effect algebras} that were
introduced by Foulis and Bennett \cite{FoBe}.  This is an algebraic
structure with a partially defined  primary operation, $+$, that
model the join of two mutually excluding quantum events. It was
inspired by algebraic properties of the effect operators on a
Hilbert space and by POV-measures.

If an effect algebra satisfies the Riesz Decomposition Property
((RDP) for short), it is always an interval in an Abelian po-group
(= partially ordered group) with interpolation and with  strong
unit, \cite{Rav}. The theory of quantum structures contains Boolean
algebras, orthomodular posets and lattices, orthoalgebras, etc. The
monograph \cite{DvPu} can serve as a comprehensive source on effect
algebras.

{\it Effect-clans} form a very important family of effect algebras.
They are effect algebras of $[0,1]$-valued functions with effect
algebraic operations defined by points. If an effect-clan is closed
with respect to pointwise limit of a sequence of nondecreasing
functions, the structure is called an {\it effect-tribe} and
effect-tribes are examples of monotone $\sigma$-complete effect
algebras. The importance of effect-tribes is underlying by the fact
that every monotone $\sigma$-complete effect algebra with (RDP) is
an epimorphic  image of some effect-tribe with (RDP), \cite{BCD}, an
analogue of the Loomis--Sikorski Theorem.

A very important subfamily of effect algebras  is a family of  {\it
MV-algebras} introduced in \cite{Cha} that model many-valued logic.
They are algebraically equivalent to Phi-symmetric effect algebras
that were introduced in \cite{BeFo}. In the theory of effect
algebras, MV-subalgebras of effect algebras represent so-called
systems of simultaneously  commensurable events of the given effect
algebra, \cite{DvPu}.

In the present paper, we will study the effect-tribes in more
details.  We will concentrate to description of states (= analogues
of probability measures) and to representation of $\sigma$-additive
states by standard integrals on some appropriate $\sigma$-algebra of
subsets closely connected with the tribe.

In addition, we show that every state on an effect algebra can be
represented by an integral through a regular Borel probability
measure on a Borel $\sigma$-algebra defined on the state space. This
extends results from \cite{Dvu4} where this was proved for interval
effect algebras (= intervals in  po-groups). Such a representation
for states on MV-algebras was studied in \cite{Kro, Pan}.

The paper is organized as follows. The elements of the theory of
effect algebras are given in Section 2. States and a new criterion
for extremal states on effect algebras are presented in Section 3.
The effect-tribes and states on them are studied in Section 4, where
it is shown a relation to Butnariu--Klement Theorem \cite{BuKl} on
states on tribes that are MV-algebras of functions. Section 5 deals
with monotone $\sigma$-complete effect algebras and we show that
every $\sigma$-convex combination of extremal states is a
Jauch-Piron state. Finally, Section 6 generalizes representation of
states on effect algebras by integrals through  regular Borel
probability measures.

\section{Effect Algebras,  Effect-Clans and Effect-Tribes}%2

According to \cite{FoBe}, an {\it effect algebra} is a partial
algebra $E = (E;+,0,1)$ with a partially defined operation $+$ and
two constant elements $0$ and $1$  such that, for all $a,b,c \in E$,
\begin{enumerate}

\item[(i)] $a+b$ is defined in $E$ if and only if $b+a$ is defined, and in
such a case $a+b = b+a;$

%\vspace{-2mm}
\item[(ii)] $a+b$ and $(a+b)+c$ are defined if and
only if $b+c$ and $a+(b+c)$ are defined, and in such a case $(a+b)+c
= a+(b+c);$

%\vspace{-2mm}
\item[(iii)] for any $a \in E$, there exists a unique
element $a' \in E$ such that $a+a'=1;$

%\vspace{-2mm}
\item[(iv)] if $a+1$ is defined in $E$, then $a=0.$
\end{enumerate}

If we define $a \le b$ if and only if there exists an element $c \in
E$ such that $a+c = b$, then $\le$ is a partial order on the set
$E$, and we write $c:=b-a.$ It is clear that $a' = 1 - a$ for any $a
\in E.$

The partial operation $+$ denotes in fact an analogue of the  join
of two mutually excluding elements. For example, every Boolean
algebra, orthomodular lattice, orthomodular poset, orthoalgebra is
an example of effect algebras, and nowadays, effect algebras are one
of the most important category of  so-called {\it quantum
structures.} Today, a modern approach to Hilbert space quantum
mechanics, POV-measures, can be elegantly described in the frame
work of effect algebras.

For a comprehensive monograph on the theory of effect algebras, see
\cite{DvPu}.

There are two very important families of effect algebras, interval
effect algebras, i.e., those that are an interval in a partially
ordered groups (po-groups for short) and effect algebras as a system
of $[0,1]$-valued functions where the effect algebraic operations
are defined by points.

Let $G$ be an Abelian po-group. An element $u\in G,$ $u\ge 1$ is
said to be a {\it strong unit} (= order unit) if for any $g \in G,$
there is an integer $n\ge 1$ such that $g\le nu,$ and the couple
$(G,u)$ is said to be a {\it unital po-group}. The interval
$\Gamma(G,u):=[0,u]$ can be endowed with the effect algebraic
operation, $+$, that is the restriction of the group addition $+$ to
the interval $[0,u];$ then $(\Gamma(G,u);+,0,u)$ is an
effect-algebra with $a'=u-a.$

An effect algebra that is either of the form $\Gamma(G,u)$ for some
element $u\ge 0$ or is isomorphic with some $\Gamma(G,u)$  is called
an {\it interval effect algebra}.

For example,  let $\mathcal B(H)$ be the system of all Hermitian
operators of a Hilbert space $H$ (real, complex or quaternionic).
Then $\mathcal B(H)$ is a po-group where the ordering of Hermitian
operators is defined by  $A \le B$ iff $(A\phi,\phi) \le
(B\phi,\phi)$ for all $\phi \in H,$ and the identity operator $I$ is
its strong unit. Then $\mathcal E(H):=\Gamma({\mathcal B}(H),I),$
the system of effect operators, is an important  example of interval
effect algebras used also in quantum mechanics, and
$$
\mathcal E(H)=\Gamma(\mathcal B(H),I).\eqno(2.1)
$$

An element $u\in G^+:=\{g\in G: g\ge 0\}$  is said to be {\it
generative} if every element $g \in G^+$ is a group sum of finitely
many elements of $\Gamma(G,u),$ and $G = G^+-G^+;$ then $u$ is a
strong unit. If $E$ is an interval effect algebra, then there is a
po-group $G$ with a generative strong unit $u$ such that $E \cong
\Gamma(G,u)$ and every $H$-valued measure $p:\Gamma(G,u) \to H$ can
be extended to a group-homomorphism $\phi$ from $G$ into $H.$ If $H$
is a po-group, then $\phi$ is a po-group-homomorphism. Then $\phi$
is unique and $(G,u)$ is also unique up to isomorphism of unital
(po-) groups, see \cite[Cor 1.4.21]{DvPu}. In such a case, the
element $u$ is said to be a {\it universal strong unit} for
$\Gamma(G,u)$ and the couple $(G,u)$ is said to be a {\it unigroup}.
For example, the identity operator $I$ is a universal strong unit
for $\Gamma({\mathcal B}(H),I),$ \cite[Cor 1.4.25]{DvPu}, similarly
for $\Gamma(A,I),$ where $A$ is a von Neumann algebra on $H.$

We recall that if an effect algebra $E$ has (RDP) if  $x_1 + x_2 =
y_1 + y_2$ implies there exist four elements $c_{11}, c_{12},
c_{21}, c_{22} \in E$ such that $x_1 = c_{11} + c_{12},$ $x_2 =
c_{21} + c_{22},$ $y_1 = c_{11} + c_{21},$ and $y_2 = c_{12} +
c_{22}.$

A partially ordered Abelian group $(G;+,0)$ (po-group in short) is
said to satisfy  {\it interpolation} provided  given
$x_1,x_2,y_1,y_2$ in $G$ such that $x_1,x_2\leq y_1,y_2$  there
exists $z$ in $G$ such that $x_1,x_2\leq z\leq y_1,y_2,$  and $G$ is
called an {\it interpolation group.}

Ravindran \cite{Rav} (\cite[Thm 1.7.17]{DvPu}) showed that every
effect algebra with (RDP) is of the form (2.1) for some
interpolation unital po-group $(G,u)$.  In other words, every effect
algebra with (RDP) is an interval effect algebra.

Bennett and Foulis, \cite{BeFo}, introduced an important subfamily
of effect algebras, Phi-symmetric effect algebras. They are
equivalent to MV-algebras. We recall that an {\it MV-algebra} is an
algebra $(M;\oplus,^*,0)$ of signature $\langle 2,1,0\rangle,$ where
$(M;\oplus,0)$ is a commutative monoid with neutral element $0$, and
for all $x,y \in M$
\begin{enumerate}
\item[(i)]  $(x^*)^*=x,$
\item[(ii)] $x\oplus 1 = 1,$ where $1=0^*,$
\item[(iii)] $x\oplus (x\oplus y^*)^* = y\oplus (y\oplus x^*)^*.$
\end{enumerate}

%We define also two  additional total operations $\odot$ and
%$\ominus$ on $M$ via $x\odot y:= (x^*\oplus y^*)^*$ and $x\ominus y
%= x\odot y^*.$

If on an MV-algebra $M$ we define a partial operation, $+$, by $a+b$
is defined in $M$ iff $a\le b^*,$ and we set then $a+ b:= a\oplus b.$
Then $(M;+,0,1)$ is an interval effect algebra with (RDP), moreover,
thanks to \cite{Mun},  every MV-algebra is in fact an interval
$\Gamma(G,u),$ where $G$ is a unital $\ell$-group (= lattice ordered
group) with strong unit $u$ with $a^*=u-a$ and $a\oplus b:=
(a+b)\wedge u,$ $a,b \in \Gamma(G,u).$

The importance of MV-algebra follows from the fact if an effect
algebra is a lattice, then it can be covered by blocks, maximal sets
of mutually compatible elements, and each of blocks is an
MV-algebra, see \cite{Rie}.

Finally we present effect algebras of $[0,1]$-valued functions with
effect algebraic operations defined by points.

Let ${\mathcal E}$ be a system of $[0,1]$-valued functions on $X\ne
\emptyset$ such that (i) $1 \in {\mathcal E}$, (ii) $f \in {\mathcal
E}$ implies $1-f \in {\mathcal E}$, and (iii) if $f,g \in {\mathcal
E}$ and $f(x) \le 1 -g(x)$ for any $x \in X$, then $f+g \in
{\mathcal E}$. Then ${\mathcal E}$ is an effect algebra of
functions, called an {\it effect-clan}, that is not necessarily a
Boolean algebra nor an MV-algebra.

A system $\mathcal E \subseteq [0,1]^X$ is said to be a {\it Bold
algebra} if (i) $1\in \mathcal E,$ (ii) $f \in \mathcal E$ implies
$1-f \in \mathcal E,$ and (iii) $f,g \in \mathcal E$ implies
$f\oplus g \in \mathcal E,$ where $(f\oplus
g)(x):=\min\{f(x)+g(x),1\},$ $x \in X.$ Then every Bold algebra is
an effect-clan with (RDP) that is an MV-algebra of functions with
pointwise defined MV-operations. For example, $\chi_A \oplus \chi_B=
\chi_{A\cup B}.$

An effect algebra is {\it monotone $\sigma$-complete}   provided
that for every ascending (descending) sequence $x_1 \le x_2 \le
\cdots $ ($x_1 \ge x_2 \ge \cdots$) in $E$ which is bounded above
(below) in $E$ has a supremum (infimum) in $E.$

An {\it effect-tribe} on a set $X \ne \emptyset$ is any system
${\mathcal T} \subseteq [0,1]^X$ such that (i) $1 \in {\mathcal T}$,
(ii) if $f \in {\mathcal T},$ then $1-f \in {\mathcal T}$, (iii) if
$f,g \in {\mathcal T}$, $f \le 1-g$, then $f+g \in {\mathcal T},$
and (iv) for any sequences $\{f_n\}$ of elements of ${\mathcal T}$
such that $f_n \nearrow f$ (pointwisely), then $f \in {\mathcal T},$
i.e. if $f_n(x) \nearrow f(x)$ for every $x \in X,$ then $f \in
{\mathcal T}.$ Every effect-tribe is a monotone $\sigma$-complete
effect algebra that is also an effect-clan.

We recall that a {\it tribe} on $X \ne \emptyset$ is a collection
${\mathcal T}$ of functions from $[0,1]^X$ such that (i) $1 \in
{\mathcal T}$, (ii) if $f \in {\mathcal T}$, then $1 - f \in
{\mathcal T},$ and (iii) if $\{f_n\}_n$ is a sequence from
${\mathcal T}$, then $\min \{\sum_{n=1}^\infty f_n,1 \}\in {\mathcal
T}.$  A tribe is always a $\sigma$-complete MV-algebra (with respect
to the pointwise operations and ordering).  For example,
$\min\{\sum_n\chi_{A_n}(x),1\}=\chi_{\bigcup_n A_n}(x), $ $x \in X.$

\section{States on Effect Algebra}%3

A {\it state} on an effect algebra $E$ is any mapping $s: \ E \to
[0,1]$ such that (i) $s(1) = 1$, and (ii) $s(a+b) = s(a) + s(b)$
whenever $a+b$ is defined in $E$.

We denote by ${\mathcal S}(E)$ the set of all states on $E$. It can
happen that an effect algebra is stateless  \cite[Ex 4.2.4]{DvPu}.
Fortunately, every interval effect algebra has at least one state,
see  \cite[Cor 4.4]{Goo}, in particular, this is true if $E$ has
(RDP).

A state $s$ is said to be {\it extremal} if $s = \lambda s_1 +
(1-\lambda)s_2$ for $\lambda \in (0,1)$ implies $s = s_1 = s_2.$ By
$\partial_e{\mathcal S}(E)$ we denote the set of all extremal states
of ${\mathcal S}(E)$ on $E.$ We say that a net of states,
$\{s_\alpha\}$, on $E$ {\it weakly converges} to a state, $s,$ on
$E$ if $s_\alpha(a) \to s(a)$ for any $a \in E$. In this topology,
${\mathcal S}(E)$ is a compact Hausdorff topological space and every
state on $E$ lies in the weak closure of the convex hull of the
extremal states as it follows from the Krein-Mil'man Theorem,
\cite[Thm 5.17]{Goo}.

We say that a system of states, ${\mathcal S},$ on an effect algebra
$E$ is (i) {\it order determining} if $s(a)\le s(b)$ for any $s \in
{\mathcal S}$ yields $a\le b,$ and (ii) {\it separating} if, for
$a,b\in E,$ $a\ne b$, there is a state $s \in {\mathcal S}$ such
that $s(a)\ne s(b),$ or, equivalently, $s(c)=s(d)$ for any $s \in
{\mathcal S}$ entails $c=d.$

Let  $\mathcal S(E)\ne \emptyset.$ Given $a \in E,$ we define the
mapping $\hat a:\mathcal S(E) \to [0,1]$  by $\hat a(s):= s(a),$ $s
\in \mathcal S(E),$  and let  $\widehat E:=\{\hat a: a \in A\}.$
Then $\widehat E$ is an effect-clan if, e.g.,  $\mathcal S(E)$ is
order determining (see \cite[Prop 4.1]{BCD1}; in such a case, $E$
and $\hat E$ are isomorphic), or $E$ is an MV-algebra; and in this
case the natural mapping $\psi(a):= \hat a$ ($a \in E$) is a
homomorphism of $E$ onto $\widehat E.$

Nevertheless that $\widehat{a+b}=\hat a+\hat b,$ $\widehat E$ is not
necessarily an effect-clan:

\begin{example}\label{ex:4.1}
Let $\mathbb Q$ be the set of all rational numbers and let
$G=\mathbb Q\times \mathbb Q$ be ordered by the strict ordering,
i.e. $(g_1,g_2)\le (h_1,h_2)$ iff $g_1 < h_1$ and $g_2< h_2$ or
$g_1=h_1$ and $g_2=h_2.$ If we set $u=(1,1),$ $(G,u)$ is a unital
po-group with interpolation.  Then for the interval effect algebra
$E=\Gamma(G,u),$ $\partial_e \mathcal S(E)=\{s_0,s_1\},$ where
$s_0(g,h):=h$ and $s_1(g,h):=g,$  $\mathcal S(E)$ is separating but
not order determining, and $\widehat E$ is not an effect-clan
because if  $a=(0.3,0.3),$ $b=(0.7,0.4),$ then $\hat a \le 1-\hat b$
but $\hat a + \hat b \notin \widehat E.$

\end{example}

Let $(G,u)$ be an Abelian po-group with strong unit. By a {\it
state} on $(G,u)$ we mean any mapping $s:\ G \to \mathbb R$ such
that (i) $s(g+h) = s(g)+ s(h)$ for all $g,h \in G,$ (ii) $s(G^+)
\subseteq \mathbb R^+$, and (iii) $s(u) = 1.$ In other words, a
state on $(G,u)$ is any po-group homomorphism from $(G,u)$ into the
po-group $(\mathbb R,1)$ that preserves fixed strong units. We
denote by ${\mathcal S}(G,u)$ the set of all states and by
$\partial_e {\mathcal S}(G,u)$ the set  of all extremal states,
respectively,  on $(G,u)$. According to \cite[Cor. 4.4]{Goo},
${\mathcal S}(G,u) \ne \emptyset$ whenever $u >0.$ In a similar way
as for effect algebras, we define the weak convergence of states on
$(G,u)$, and analogously, ${\mathcal S}(G,u)$ is a compact Hausdorff
topological space and every state on $(G,u)$ is a weak limit of a
net of convex combinations from $\partial_e {\mathcal S}(G,u).$

If $E= \Gamma(G,u)$, where $(G,u)$ is a unigroup, then  every state
on $E$ can be extended to a unique state on $(G,u)$, and conversely,
the restriction of any state on $(G,u)$ to $E$ gives a state on $E$.
In addition, extremal states on $E$ are the restrictions of extremal
states on $(G,u)$,  the space ${\mathcal S}(E)$ is affinely
homeomorphic with ${\mathcal S}(G,u)$, and the space $\partial_e
{\mathcal S}(E)$ is homeomorphic with $\partial_e {\mathcal
S}(G,u)$. We recall that a mapping from one convex set into another
convex set is {\it affine} if it preserves all convex combinations.

We recall that if $E$ is an MV-algebra, then a state $s$ on $E$ is
extremal if and only if $s(x\wedge y) =\min\{s(x),s(y)\}$, $x,y \in
E$, therefore $\partial_e {\mathcal S}(E)$ is compact (see, e.g.
\cite[Prop. 6.1.19]{DvPu},  \cite[Thm 12.18]{Goo}).

On the other hand,  if $E$ is an effect algebra with (RDP), then
$\partial_e {\mathcal S}(E)$ is not necessarily compact,   see for
example \cite[Exam. 6.10]{Goo}.  We note that due to \cite[Thm
12.14]{Goo}, a state $s$ on an interpolation unital po-group $(G,u)$
is extremal iff

$$\min\{s(x),s(y)\}=\sup\{s(z):\ z \in G^+,\ z\le x,z\le y\} \eqno(3.1)$$
for all $x,y \in G^+.$

We recall that if $s$ is an extremal state on $E$ with (RDP), then
$$\min\{s(x),s(y)\}=\sup\{s(z):\ z \in E,\ z\le x,z\le y\} \eqno(3.2)$$
for all $x,y \in E.$

We note that we do not know whether (3.2) implies (3.1). Partial
positive answers are given in Proposition \ref{pr:3.2} and
Proposition \ref{pr:8.7.1} below.

We say   an effect algebra $E$ is {\it divisible} if for every
positive integer $n$ and $x \in E,$ there is an element $y \in E$
such that $ny =x$ and we write $y = \frac{1}{n}x.$ A po-group $G$ is
(i) {\it divisible} if for every positive integer $n$ and $g \in G,$
there is an element $h \in G$ such that $nh =g,$ (ii) {\it
unperforated} if $ng \ge 0$ for some integer $n\ge 0$ implies $x \ge
0.$

\begin{proposition}\label{pr:3.1} If an effect algebra $E$ with {\rm
(RDP)} is  divisible, then a  unital po-group $(G,u)$ with
interpolation such that $E=\Gamma(G,u)$ is divisible.  If an
unperforated unital po-group $(G,u)$ with interpolation is
divisible, then $E=\Gamma(G,u)$ is divisible.
\end{proposition}

\begin{proof} Let $E$ be divisible. Given $g \in G^+,$ there is an integer
$k\ge 1$ such that $g \le ku.$  Then $g = g_1+\cdots +g_k$ where
each $g_i \le u.$ For every $n\ge 1,$ $\frac{1}{n}g_i \in E,$ so
that $g = n \frac{1}{n}g_1+\cdots + n \frac{1}{n} g_k =
n(\frac{1}{n}g_1+\cdots +\frac{1}{n}g_k)$.  Finally, if $g \in G$ is
arbitrary, then $g = g_1-g_2,$ where $g_1,g_2 \ge 0.$  Then
$\frac{1}{n}g_1, \frac{1}{n}g_2 \in G$ and hence $g= ng',$ where $g'
= \frac{1}{n}g_1- \frac{1}{n}g_2 \in G.$

Now if $G$ is divisible, then for every $x\in E=\Gamma(G,u)$ and
$n\ge 1$ there is an element $y\in G$ such that $ny=x\ge 0,$ hence
$y\ge 0$ so that $y \in E$ and $E$ is divisible.
\end{proof}

The following result is true also for $E=\Gamma(G,u),$ where $(G,u)$
is a unigroup.

\begin{proposition}\label{pr:3.2} Let $E$ be a divisible effect
algebra with {\rm (RDP)}.  Then a state $s$ on $E$ is extremal if
and only if {\rm (3.2)} holds.
\end{proposition}

\begin{proof}
Let $E=\Gamma(G,u)$ for some unital po-group $(G,u)$ with
interpolation. If $s$ is extremal then its unique restriction to
$(G,u)$ is again extremal so that we have  (3.1) and consequently
(3.2).

Conversely, let (3.2) hold and assume that $s$ is not extremal. Then
$s = \lambda_1 s_1 + \lambda_2s_2$ for some distinct states $s_1$
and $s_2$ and real numbers $\lambda_1,\lambda_2 >0$ such that
$\lambda_1+\lambda_2 =1.$  We will write $s,s_1,s_2$ also for their
unique extension on $(G,u).$

Therefore, there is an element $a \in E$ such that $s_1(a) \ne
s_2(a).$ Without loss of generality, we can assume $0\le s_1(a) <
s_2(a).$ There are two positive integers $m$ and $n$ such that
$s_1(a) < m/n <s_2(a).$  Put $x = na$ and $y = mu.$  Then $s_1(x) <
s_1(y)$ and $s_2(y) < s_2(x).$  Hence $\frac{1}{nm} x = \frac{1}{m}a
\le a$ and $\frac{1}{nm}y \le \frac{1}{n}u \le u.$ For any $z\in E$
with $z \le \frac{1}{m}a$ and $z\le \frac{1}{n}u,$ we have
$$s(z) =
\lambda_1 s_1(z) + \lambda_2 s_2(z) \le \lambda_1 s_1(\tfrac{1}{m}a)
+ \lambda_2 s_2(\tfrac{1}{n}u).
$$
Hence,
$$\min\{s(\tfrac{1}{m}a), s(\tfrac{1}{n}u)\} = \sup\{s(z): z
\le \tfrac{1}{m}a, \ z\le \tfrac{1}{n}u\} \le \lambda_1
s_1(\tfrac{1}{m}a) + \lambda_2 s_2(\tfrac{1}{n}u).
$$
Since $\lambda_1
>0$ and $s_1(x)< s_1(y),$ we have

\begin{eqnarray*}
\lambda_1 s_1(\tfrac{1}{m}a) + \lambda_2 s_2(\tfrac{1}{n}u) &=&
\lambda_1 s_1(\tfrac{1}{nm} x) + \lambda_2 s_2(\tfrac{1}{nm}
y)\\
&=& \tfrac{1}{nm}(\lambda_1 s_1(x) + \lambda_2 s_2(y))\\
&<& \tfrac{1}{nm}(\lambda_1 s_1(y) + \lambda_2 s_2(y)) \\
&=& \tfrac{1}{nm}s(y) = s(\tfrac{1}{n}u).
\end{eqnarray*}

Similarly, $\lambda_1 s_1(\tfrac{1}{m}a) + \lambda_2
s_2(\tfrac{1}{n}u)< s(\tfrac{1}{m}a)$ which gives a contradiction.
 \end{proof}

It is important to recall that according to a delicate result of
Choquet \cite[Thm I.5.13]{Alf}, $\partial_e {\mathcal S}(E)$ is
always a Baire space in the relative topology induced by the
topology of ${\mathcal S}(E)$, i.e. the Baire Category Theorem holds
for $\partial_e {\mathcal S}(E).$

A $\sigma$-additive state on an effect algebra $E$ is a state $s$
such that if $a_n \le a_{n+1}$ for any $n\ge 1$  and $a=\bigvee_n
a_n \in E$ (and we write $a_n \nearrow a$), then $s(a) = \lim_n
s(a_n).$ It is easy to verify that a state $s$ is $\sigma$-additive
iff $a_n \searrow 0$ then $\lim_n s(a_n)=0.$

Let $H$ be a separable Hilbert space (real, complex or quaternionic)
with an inner product $(\cdot,\cdot)$, and ${\mathcal L}(H)$ be the
system of all closed subspaces of $H.$ Then ${\mathcal L}(H)$ is a
complete orthomodular lattice \cite{Dvu}. Given a unit vector $x \in
H$, let
$$
p_x(M):= (P_Mx,x),\quad M \in {\mathcal L}(H),
$$
where $P_M$ is the orthogonal projector of $H$ onto $M.$ Then $p_x$
is a $\sigma$-additive state on ${\mathcal L}(H),$ called a {\it
pure state}. The system of all pure states is order determining.  If
$T$ is a positive Hermitian operator of finite trace (i.e.
$\sum_i(Tx_i,x_i) <\infty$ for any orthonormal basis $\{x_i\}$ of
$H$, and we define $\mbox{tr}(T) :=\sum_i(Tx_i,x_i)$) such that
tr$(T)=1,$ then
$$
s_T(M):=\mbox{tr}(TP_M),\quad M \in {\mathcal L}(H),\eqno(3.3)
$$
is a $\sigma$-additive state, and according to Gleason's theorem,
\cite{Gle}, \cite[Thm 3.2.24]{Dvu}, if $\dim H\ge 3,$ for every
$\sigma$-additive state $s$ on ${\mathcal L}(H)$, there is a unique
positive Hermitian operator $T$ with tr$(T)=1$ such that $s=s_T.$

If $s$ is a state on $\mathcal L(H),$ $3\le\dim H \le \aleph_0$,
then the state $s$ can be uniquely extended to a state, $\hat s$ on
$\mathcal B(H),$ moreover, $\hat s(\alpha A)=\alpha \hat s(A)$ for
any $\alpha \in \mathbb C$ ($\alpha \in \mathbb R$) and $A \in
\mathcal B(H).$

If $\dim H = 2,$ there are two-valued states on $\mathcal L(H)$ and
they are extremal.

We denote by $\mathcal S_c(\mathcal L(H))$ the system of all states
that can be extended to a state on $\mathcal B(H).$  This is also a
convex state, and if $\dim H \ge 3,$ $\mathcal S_c(\mathcal L(H)) =
\mathcal S(\mathcal L(H)).$

Now if $s$ is a finitely additive state on ${\mathcal L}(H),$ $\dim
H \ge 3,$ then due to the Aarnes theorem,  \cite[Thm 3.2.28]{Dvu},
$s$ can be uniquely expressed in the form
$$ s = \lambda s_1 +(1-\lambda)s_2,$$
where $s_1$ is a $\sigma$-additive state and $s_2$ is a finitely
additive state vanishing on each finite-dimensional subspace of $H.$

Moreover, a pure state $p_x$ is an extreme point of the set of
$\sigma$-additive states, as well as it is an extremal state of
${\mathcal L}(H).$ Other extremal states of ${\mathcal L}(H)$ vanish
on each finite-dimensional subspace of $H.$  Since every state on
${\mathcal L}(H)$ can be  extended into a unique state on ${\mathcal
B}(H),$ see e.g. \cite[Thm 3.3.1]{Dvu},  the state spaces ${\mathcal
S}({\mathcal L}(H)),$ ${\mathcal S}({\mathcal E}(H))$ and ${\mathcal
S}({\mathcal B}(H))$ are mutually affinely homeomorphic whenever
$\dim H \ge 3.$

Let us set $\Omega(H):=\{x \in H: ||x||= 1\}$ and for any $A \in
\mathcal E(H)$ we define $f_A(x):=(Ax,x),$ $x \in \Omega(H).$  Then
$\mathcal T(H):=\{f_A: A \in \mathcal E(H)\}$ is an effect-tribe.

In addition, the set $\mathcal S=\{p_x: x \in \Omega(H)\}$ is an
order determining system of states. The effect-algebra $\mathcal
E(H)$ is monotone $\sigma$-complete and it is isomorphic with the
effect-tribe $\widehat {\mathcal E(H)}$ and as well as with
$\mathcal T(H).$

\section{States on Effect-Tribes as Integrals}%4

We will characterize $\sigma$-states on effect-tribes by integrals
through probability measures on a special $\sigma$-algebra of sets
connected with the given tribe. We show how we can generalize the
Butnariu--Klement Theorem, \cite{BuKl}, that was proved for tribes.

Let $E$ be an effect algebra. For a given  element $e \in E$, we
denote by $[0,e] := \{x\in E:\, 0\le x \le e\}.$ Then $[0,e]$
endowed with $+$ restricted to $[0,e]\times [0,e]$ is an effect
algebra $[0,e] =([0,e];+,0,e)$. For any $x \in [0,e]$ we have
$x^{'_e} := e -x.$

An element $e$ of an effect algebra $E$ is said to be {\it central}
(or {\it Boolean}) if there exists an isomorphism
$$
f_e:\, E \to [0,e] \times [0,e'] \eqno(4.1)
$$
such that $f_e(e) =(e,0)$ and if $f_e(x) = (x_1,x_2)$, then $x= x_1
+ x_2$ for any $x \in E.$

We denote by $B(E)$ the set of all central elements of $E$, and
$B(E)$ is said to be the {\it center} of $E$. We recall that $0,1
\in B(E).$ According to \cite[Thm 1.9.14]{DvPu}, $B(E)$ is an effect
subalgebra of $E$ and it is a Boolean algebra. For any $x\in E$ and
any $e \in B(E)$ we have
$$x = (x\wedge e)+(x\wedge e'); \eqno(4.2)
$$
we define $p_e:\ E\to [0,e]$ by $p_e(x) = x\wedge e.$

If $E$ satisfies (RDP), then an element $e$ is central iff $e\wedge
e' = 0,$  \cite[Thm 3.2]{Dvu1}. We note that without (RDP), the last
statement is not always true: Take $X=\{1,...,6\}$ and let
${\mathcal T}$ be an effect-clan that  consists of the
characteristic functions of all subsets of $X$ with even numbers of
elements. Then ${\mathcal T}$ is without (RDP), any element of
${\mathcal T}$ is a characteristic function but  it is not always
central, e.g. $\chi_{\{1,2,3,4\}}\wedge \chi_{\{1,2,3,5\}} \notin
{\mathcal T}.$ Moreover, $B(\mathcal T) = \{0,1\}.$

In addition, for an effect-clan ${\mathcal T}$  with (RDP),  a
characteristic function $e\in {\mathcal T}$ is always a central
element of $\mathcal T.$ If $\mathcal T$ is even a Bold algebra,
$e\in \mathcal T$ is central iff $e$ is a characteristic function of
some subset of $X.$ As we have already seen just after (4.2),
without (RDP), the last statement is not always true.

In addition, if $E$ is monotone $\sigma$-complete, then $B(E)$ is a
Boolean $\sigma$-algebra, \cite[Thm 5.11]{Dvu1}, and for any $x\in
E$ and $\{e_n\}$ from $B(E)$ we have $x \wedge \bigvee_{n=1}^\infty
e_n= \bigvee_{n=1}^\infty (x\wedge e_n).$

Let ${\mathcal B_0}({\mathcal T})$ be the system of all subsets $A$
of $X$ such that $\chi_A$ is a central element of the effect-clan
${\mathcal T},$ i.e.
$$
{\mathcal B_0}({\mathcal T}) =\{A\subseteq X:\ \chi_A \in
B({\mathcal T})\}.
$$

\begin{lemma}\label{le:9.1'}  Let $\mathcal T$ be an effect-clan.
If $f \in B(\mathcal T)$ and $g\in \mathcal T,$ then $gf \in
\mathcal T$ and $g\wedge f= gf,$ where $gf$ denotes  the product of
two functions $g$ and $f.$

\end{lemma}

\begin{proof} Let $f=\chi_A\in B(\mathcal T).$ We have $(g\wedge f)(x) \le
g(x)=g(x)f(x)$ if $x \in A,$ and $(g\wedge f)(x) \le
f(x)=g(x)f(x)=0$ if $x \notin A.$ Similarly we have also $(g\wedge
(1-f))(x)\le g(x)(1-f)(x)$ for any $x \in X.$

Hence, $g(x)= (g\wedge f)(x) + (g\wedge (1-f))(x) \le
g(x)f(x)+g(x)(1-f)(x) = g(x),$ $x \in X,$ which proves $g\wedge f =
gf,$ and $gf \in \mathcal T.$
\end{proof}

The system ${\mathcal B_0}({\mathcal T})$ is a $\sigma$-algebra of
subsets of $X$ whenever $\mathcal T$ is an effect-tribe:

\begin{proposition}\label{pr:8.1}  Let  ${\mathcal T}$ be an
effect-clan.  Then  ${\mathcal B_0}({\mathcal T})$ is an algebra of
subsets, and if $E, F \in {\mathcal B_0}({\mathcal T}),$ then
$\chi_E\wedge \chi_F = \chi_{E\cap F}$ and $\chi_E\vee \chi_F =
\chi_{E\cup F}.$

If, in addition, $\mathcal T$ is an effect-tribe, then ${\mathcal
B}({\mathcal T})$ is a $\sigma$-algebra of subsets, and if $\{E_n\}$
is a sequence of elements from ${\mathcal B_0}({\mathcal T})$, then
$\bigwedge_n \chi_{E_n} = \chi_{\bigcap_n E_n}$ and $\bigvee_n
\chi_{E_n} = \chi_{\bigcup_n E_n}.$
\end{proposition}

\begin{proof}  Since $\chi_E$ and $\chi_F$ are central elements, by
Lemma \ref{le:9.1'}, $\chi_E \wedge \chi_F = \chi_E  \chi_F =
\chi_{E\cap F}.$ Passing to negations, we have the second identity.

Now let $\mathcal T$ be an effect-tribe and choose a sequence
$\{E_n\}$ from ${\mathcal B_0}({\mathcal T}).$  According to the
first part, we can assume that $E_n \subseteq E_{n+1}$ and let $E=
\bigcup_n E_n.$ Then $h = \bigvee_n \chi_{E_n} \in {\mathcal T}$
and, for any $x \in X,$ we have $h(x) = \lim_n \chi_{E_n}(x) =
\chi_E(x).$  Similarly we can prove by Lemma \ref{le:9.1'} that
$f\chi_E = \lim_n f\chi_{E_n} \in \mathcal T$ which proves $f\chi_E
\le f$ and $f\chi_e\le \chi_E.$  Now if $g\le f,\chi_E$ for some $g
\in \mathcal T,$ then $g\le f\chi_E$ and this yields $f\wedge \chi_E
=f\chi_E.$  In a similar way we can prove that $f\wedge
(1-\chi_E)=f(1-\chi_E)\in \mathcal T$  so that $f = (f\wedge
\chi_E)+(f\wedge (1-\chi_E))$ and hence $\chi_E\in B(\mathcal T)$
and $E \in {\mathcal B_0}({\mathcal T}).$

The second equality is now evident.
\end{proof}

It is worthy to recall that if $\mathcal T$ is an effect-clan with
(RDP) then not every $f \in B({\mathcal T)}$ is a characteristic
function. In fact, let $X =[0,1]$ and ${\mathcal
T}=\{0,\mbox{id}_{[0,1]}, 1- \mbox{id}_{[0,1]}, 1\}.$  Then
${\mathcal T}$ is an MV-algebra (in fact a Boolean algebra) such
that $B(\mathcal T) = \mathcal T$ and ${\mathcal B_0}(\mathcal
T)=\{0,1\}.$

At any rate, ${\mathcal B_0}(\mathcal T)$ is a Boolean
sub-$\sigma$-algebra of $B(\mathcal T)$ whenever $\mathcal T$ is an
effect-tribe.

Let $\mathcal T$ be an effect-clan of functions on $X\ne \emptyset,$
then
$$\mathcal S_0(\mathcal T):=\{A\subseteq X: \chi_A \in \mathcal T\}
$$
is a system of subsets such that  {\rm (i)} $\emptyset, X\in
\mathcal S_0(\mathcal T),$ {\rm (ii)} if $A \in \mathcal
S_0(\mathcal T),$ then $X \setminus A \in \mathcal S_0(\mathcal T),$
and, in addition if $\mathcal T$ is an effect-tribe {\rm (iii)} if
$A_n\in \mathcal S_0(\mathcal T)$, $A_n\subseteq A_{n+1}$ for any
$n\ge 1,$ then $A = \bigcup_n A_n \in \mathcal S_0(\mathcal T).$ The
system $\mathcal S_0(\mathcal T)$ is not necessarily an algebra but
$$
\mathcal B_0\mathcal {(T)} \subseteq \mathcal S_0(\mathcal T).
$$

If $\mathcal T$ satisfies {\rm (RDP)}, then $\mathcal S_0(\mathcal
T)$ is an algebra, because then any element $a=\chi_A \in \mathcal
T$ is central, hence, $A \in \mathcal S_0(\mathcal T),$ and
$$
\mathcal B_0\mathcal {(T)} = \mathcal S_0(\mathcal T).
$$
Hence, $\mathcal S_0(\mathcal T)$ is a $\sigma$-algebra of subsets
whenever $\mathcal T$ is an effect-tribe. Indeed, because $\mathcal
S_0(\mathcal T)$ is a Boolean algebra, without loss of generality,
we can assume that $\{A_n\}$ is a nondecreasing. Then there is $a
\in \mathcal T$ such that $\chi_{A_n}\nearrow a$ by points.  Hence,
$a(x)\in \{0,1\}$ for any $x \in X,$ so that $a=\chi_A$ for some $A
\subseteq X$ and $\bigcup_n A_n = A \in \mathcal S_0(\mathcal T).$

If $s$ is a  finitely additive state or a  $\sigma$-additive state
on ${\mathcal T},$ then the mapping $\mu_s$ defined on ${\mathcal
B}(\mathcal T)$ by
$$\mu_s(A):=
s(\chi_A),\ A \in {\mathcal B_0}({\mathcal T}),\eqno(4.3)
$$
is a  finitely additive normalized measure, or a $\sigma$-additive
one on ${\mathcal B_0}({\mathcal T}).$

Now we define  a  special type of effect-tribes.  We recall that a
{\it tribe} on $X \ne \emptyset$ is a collection ${\mathcal T}$ of
functions from $[0,1]^X$ such that (i) $1 \in {\mathcal T}$, (ii) if
$f \in {\mathcal T}$, then $1 - f \in {\mathcal T},$ and (iii) if
$\{f_n\}_n$ is a sequence from ${\mathcal T}$, then $\min
\{\sum_{n=1}^\infty f_n,1 \}\in {\mathcal T}.$  A tribe is always a
$\sigma$-complete MV-algebra (with respect to the pointwise
operations and ordering).

Butnariu and Klement, \cite{BuKl}, showed that if ${\mathcal T}$ is
a tribe, then every element $f \in {\mathcal T}$ is measurable with
respect to ${\mathcal B_0}({\mathcal T})$ and if $s$ is a
$\sigma$-additive state on ${\mathcal T},$ then

$$s(f) = \int_X f \dx \mu_s, \ f \in {\mathcal T}, \eqno(4.4)
$$
where $\mu_s$ is as in (4.3).

It is worthy to note that  in the case of the tribe ${\mathcal
T(H)}$ we have $ B({\mathcal T(H)})=\{0,1\}$ and $\mathcal
B_0({\mathcal T(H)})$ is trivial, so that only constant functions
from ${\mathcal T(H)}$ are $\mathcal B_0({\mathcal
T(H)})$-measurable, so that the Butnariu--Klement result is not true
for all effect-tribes.

In what follows, we show that there is even an effect-tribe
$\mathcal T$ with (RDP) that is  not a lattice and not every element
of $\mathcal T$ is measurable with respect to ${\mathcal
B_0}({\mathcal T}),$ so that the Butnariu--Klement theorem is not
necessarily valid for all effect-tribes with (RDP):

\begin{example}\label{ex:8.2}
There is an effect-tribe ${\mathcal T}$  with {\rm (RDP)} over a set
$X\ne \emptyset$ such  that $\mathcal T$ is not a lattice, and not
every element is measurable with respect to ${\mathcal
B_0}({\mathcal T}).$ Moreover, let $\mathcal S_0$ be the least
$\sigma$-algebra of subsets of $X$ such that each $f \in \mathcal T$
is $\mathcal S_0$-measurable. There are two different probability
measures, $\mu_0$ and $\mu_1,$ on $\mathcal S_0$ such that formula

$$
s_\mu(f):=\int_X f \dx \mu,\ f \in {\mathcal T},\eqno(4.5)
$$
defines the same $\sigma$-additive state on ${\mathcal T}$ if we set
$\mu=\mu_0$ and  $\mu=\mu_1.$

\end{example}

\begin{proof}  Let $X$ be an uncountable set with two distinct  elements
$a,b \in X.$  Let $G$ be the set of those bounded functions $f:\ X
\to \mathbb R$ satisfying $f(x)= (f(a)+f(b))/2$ for all but
countably many $x \in X.$ If $1$ is a constant function $1$, then
$1$ is a strong unit in $G,$ and according to \cite[Ex. 16.1, Ex.
16.8]{Goo}, $G$ is an interpolation group, monotone
$\sigma$-complete but not lattice-ordered.

Then ${\mathcal T}:= \Gamma(G,1)$ is with (RDP), where $1$ is the
function equals $1.$ If $\{f_n\}$ is a monotone sequence from
$\mathcal T,$ then $f(x)=\lim_n f_n(x),$ $x \in X,$ is a function
from $\mathcal T,$ so that $\mathcal T$ is an effect-tribe.

Given $x \in X,$ $s_x(f) :=f(x)$ for every $f \in {\mathcal T},$ is
a $\sigma$-additive state on ${\mathcal T}.$

Let ${\mathcal S}$ be the $\sigma$-algebra of all subsets $A$ of $X$
such that either $A$ is countable or $X\setminus A$ is countable,
and let $\mu_0$ be a two-valued mapping on ${\mathcal S}$ such that
$\mu_0(A)=0$ iff $A$ is countable, otherwise $\mu_0(A)=1;$ then
$\mu_0$ is an extremal probability measure on ${\mathcal S}.$

Denote by ${\mathcal S}_{a,b}$ the set of all subsets $A \in
{\mathcal S}$ such that if $A$ is countable then $a,b \not\in A$ and
if $A$ is uncountable then $a,b \in A;$  ${\mathcal S}_{a,b}$ is a
$\sigma$-algebra.

In what follows, we show that  $f \in B(\mathcal T)$ iff $f =
\chi_A$ where $A \in \mathcal S_{a,b},$ and $ \mathcal S_{a,b}
={\mathcal B_0}({\mathcal T}).$

Given $f \in {\mathcal T},$ let $A_f =\{x\in X:\
f(x)=(f(a)+f(b))/2\}$ and given $\alpha \in [0,1],$ let
$A_f(\alpha)=\{x\in X:\ f(x) \le \alpha\}.$  Then $A_f = A_{1-f}.$
Moreover, $a\in A_f$ iff $b \in A_f$ iff $f(a)=f(b).$

Suppose $f \in B(\mathcal T)$. Since $\mathcal T$ has (RDP), we have
that $f \wedge (1-f)=0.$ The functions $f$ and $1-f$ are constant on
$A_f.$ Now let  $x_0 \in A_f.$ If $0<f(x_0)<1$ define $h_{x_0}:X \to
\mathbb R $ such that $h_{x_0}(x)= \min\{f(x),1-f(x)\}/2 >0$ for $x
\in A_f$ and $h_{x_0}(x) =0$ for $x \notin A_f$. Then $h_{x_0}\in
{\mathcal T}.$ Therefore, $h_{x_0}\le f,1-f$ and this contradicts
the fact that $f\wedge (1-f)=0.$ Consequently, either
$f\upharpoonright A_f =0$ or $(1-f) \upharpoonright  A_f=0.$ Now let
$x_0 \notin A_f$ then again we can show that $f(x_0) \in \{0,1\},$
so that $f$ is a characteristic function of a set $A \in {\mathcal
B}(\mathcal T)$ and if $a,b \in A,$ then $A=A_f$ and if $a,b \notin
A$, then $X\setminus A = A_f$ which yields $A \in {\mathcal
S}_{a,b}.$

Conversely, if $A \in {\mathcal S}_{a,b}$ then $f=\chi_A \in
\mathcal T$ and, therefore, $f \in B(\mathcal T)$ and $A \in
{\mathcal B_0}(\mathcal T).$

In particular, we have proved $ \mathcal S_{a,b} ={\mathcal
B}({\mathcal T}).$

From this we have (i) every $f \in {\mathcal T}$ is ${\mathcal
S}$-measurable because $A_f(\alpha)\in \mathcal S$ for any
$\alpha\in [0,1].$  (ii) Not every every $f \in {\mathcal T}$ is
$\mathcal B_0\mathcal {(T)}$-measurable: Indeed, choose $f\in
{\mathcal T}$ such that $f(a)\ne f(b).$ Then $a,b \not\in A_f,$ but
$A_f$ is uncountable so that $A_f \not\in {\mathcal S}_{a,b}.$ (iii)
$\mathcal S=\mathcal S_0.$

If $\mu$ is an arbitrary probability measure on ${\mathcal S}$, then
$$
s_\mu(f):=\int_X f \dx \mu,\ f \in {\mathcal T},
$$
is a $\sigma$-additive state on ${\mathcal T}.$ Of course, $s_x =
s_{\delta_x}$ for any $x \in X,$ where $\delta_x$ a Dirac measure
concentrated at the point $x,$ i.e. $\delta_x(A)=1$ iff $x \in A$
and $\delta_x(A)=0$ iff $x\notin A.$ In particular, $s_{\mu_0}(f)=
(f(a)+f(b))/2,$ $f \in {\mathcal T},$ so that $s_{\mu_0} =
(s_a+s_b)/2$ and $s_{\mu_0}$ is not an extremal state on ${\mathcal
T}$ however $\mu_0$ does on ${\mathcal S}.$ In addition,
$$ s_{\mu_0}(f) =\int_X f \dx \mu_0=\int_X f \dx ((\delta_a + \delta_b)/2),
\ f \in {\mathcal T},\eqno(4.6)
$$
so that the state $s_{\mu_0}$ has two different representations via
(4.5): by $\mu_0$ and by $\mu_1=(\delta_a + \delta_b)/2.$
%Since any (extremal) state $s$ on ${\mathcal T}=\Gamma(G,1)$ can be
%uniquely extended to an (extremal) state on the interpolation unital
%group $(G,1),$ to show that $s_x$ is extremal, we use the criterion
%(3.1) dealing on $(G,1).$
 \end{proof}

Now we present a generalization of a result by Butnariu and Klement
\cite{BuKl} for effect-tribes with (RDP).

\begin{theorem}\label{th:8.3}  Let $s$ be a $\sigma$-additive state
on an effect-tribe ${\mathcal T}$  with {\rm (RDP)} of functions on
a set $X\ne\emptyset$ such that every $f\in\mathcal T$ is ${\mathcal
B_0}(\mathcal T)$-measurable. Then there is a unique probability
measure $\mu$ on ${\mathcal B}(\mathcal T)$ such that
$$s(f)=\int_X f \dx \mu,\ f \in \mathcal T.\eqno(4.7)
$$
Moreover, $\mu = \mu_s,$ where $\mu_s$ is as in $(4.3).$

\end{theorem}

\begin{proof}
Set $\mu = \mu_s,$ where $\mu_s$ is defined by (4.3).  Since
${\mathcal T}$ satisfies (RDP), we can assume that ${\mathcal T}
=\Gamma(G,1)$ where $(G,1)$ is an interpolation unital po-group of
bounded functions on $X.$ The hypothesis entails that every $f\in G$
is ${\mathcal B_0}({\mathcal T})$-measurable and $s$ can be extended
to a state $\hat s$ on $(G,u).$

(1)  Assume that $f = \chi_A$ is a central element of ${\mathcal
T}.$ Then $s(f) = s(\chi_A) = \int_X \chi_A \dx \mu_s.$

(2) Let $f= \alpha \chi_A$ where $A \in {\mathcal B_0}({\mathcal
T})$ and $\alpha \in [0,1].$ First, let $\alpha = p/q$ for $p,q \in
\mathbb N.$ Then $qf = p\chi_A$ so that $\hat s(qf)=\hat s(p\chi_A)$
and $qs(f)=ps(\chi_A) = p\int_X \chi_A$ that implies $s(f)=
\frac{p}{q} s(\chi_A) = \frac{p}{q}\int_X \chi_A \dx \mu_s.$

If $\alpha$ is irrational, there are two sequence of positive
rational numbers, $\{p_n/q_n\}$ and $\{p'_n/q'_n\},$ such that $0\le
p_n/q_n \nearrow \alpha \swarrow p'_n/q'_n\le 1.$  Then
\begin{eqnarray*}
&p_nq_n' \le q_nq_n' \alpha \le q_np_n'\\
&p_nq_n'\chi_A \le q_nq_n' \alpha \chi_A \le q_np_n' \chi_A\\
&p_nq_n's(\chi_A) \le q_nq_n's(\chi_A) \alpha \le q_np_n'
s(\chi_A)\\
&\frac{p_n}{q_n}\chi_A \le  \alpha \chi_A \le \frac{p_n'}{q_n'}
\chi_A,
\end{eqnarray*}
so that $s(\alpha \chi_A) = \alpha s(\chi_A) = \int_X \alpha\chi_A
\dx \mu_s.$

(3) Let $f \in {\mathcal T}.$ For any integer $n \ge 1$ and any
$i=0,1\ldots, 2^n-1,$ we define $A_n^i=\{x\in X:\ \tfrac{i}{2^n}< f
\le \tfrac{i+1}{2^n}\} \in {\mathcal B_0}(\mathcal T)$ and
$$
f_n:=\sum_{i=0}^{2^n -1} \tfrac{i}{2_n}\chi_{A_n^i}, \quad g_n :=
\sum_{i=0}^{2^n -1} \tfrac{i+1}{2_n}\chi_{A_n^i}.
$$
Then $2^nf_n, 2^ng_n \in G,$  $ f_n \le f \le 2^ng_n,$  $2^nf_n \le
2^n f \le 2^ng_n$ and $2^nf_n, 2^nf, 2^ng_n \in G.$ Hence,
\begin{eqnarray*}
& 2^n \sum_{i=0}^{2^n-1} i\chi_{A^i_n} \le 2^nf \le 2^n
\sum_{i=0}^{2^n-1}(i+1)\chi_{A_n^i}\\
& 2^n \sum_{i=0}^{2^n-1} is(\chi_{A^i_n}) \le 2^ns(f) \le 2^n
\sum_{i=0}^{2^n-1}(i+1)s(\chi_{A_n^i})\\
&2^n\int_X f_n \dx \mu_s \le s(f) \le 2^n \int_X g_n \dx\mu_s\\
&\int_X f_n \dx\mu_s \le s(f) \le \int_X g_n \dx\mu_s.
\end{eqnarray*}
Using the Lebesgue Convergence Theorem, \cite{Hal}, we have $\int_X
\dx\mu_s \le s(f) \le \int_X f \dx\mu_s.$

Finally, let $\mu$ be any measure on ${\mathcal B_0}(\mathcal T)$
such that (4.7) holds. Then $\mu(\chi_A) = \int_X \chi_A \dx\mu =
\int_X \chi_A \dx\mu_s = \mu_s(A)$ for any $A \in {\mathcal
B_0}(\mathcal T).$
\end{proof}

\begin{corollary}\label{co:8.4} Let ${\mathcal T}$ be an
effect-tribe of functions on $X\ne \emptyset$ satisfying {\rm (RDP)}
and let $\mathcal T=\Gamma(G,1),$ where $(G,1)$ is an interpolation
unital po-group of bounded functions on a set $X\ne \emptyset   .$
Then every $\sigma$-additive state $s$ on $\mathcal T$ can be
extended to a unique $\sigma$-additive state $\hat s$ on $(G,1).$
Moreover,

$$\hat s(f)=\int_X f \dx \mu_s, \ f \in G,\eqno(4.8)
$$
where $\mu_s$ is defined by $(4.3).$
\end{corollary}

\begin{proof}  According to Theorem \ref{th:8.3}, we have (4.7).  Since
${\mathcal T}$ generates $(G,1)$ and it satisfies (RDP), the formula
(4.7) holds for $s = \hat s$ and for any $f \in G.$

The $\sigma$-additivity of $\hat s$ on $G$ means that if $ \{f_n\}$
is a monotone sequence of elements of $G$ such that $f_n \nearrow f
\in G,$ then $\lim_n \hat s(f_n) = \hat s(f).$  But this is
guaranteed by the Lebesgue Convergence Theorem.
\end{proof}

\begin{remark}\label{re:uni}{\rm  Theorem \ref{th:8.3} and Corollary
\ref{co:8.4} hold also for the case that the effect-tribe $\mathcal
T$ is an interval in some unigroup $(G,u)$ of bounded functions for
$\mathcal T$ because in this case every state on $\mathcal T$ can be
uniquely extended to a state on $(G,u).$
}
\end{remark}

%A monotone $\sigma$-complete  effect algebra with (RDP) has general
%comparability iff $\partial_e{\mathcal S}(E)$ is closed, \cite[Cor
%16.28]{Goo}, iff $E$ is a lattice iff $E$ is an MV-algebra (where
%original $+$ and derived one from $\oplus$ coincide).

\begin{proposition}\label{pr:8.5}
Let $\mathcal T$ be an effect-tribe of functions on $X\ne \emptyset$
satisfying {\rm (RDP)}. Let $\mathcal T'$ be the set of all
functions $f \in \mathcal T$ such that $f$ is $\mathcal B_0\mathcal
{(T)}$-measurable. Then $\mathcal T'$ is an effect-tribe and
$\mathcal B_0\mathcal {(T)}=\mathcal B_0\mathcal {(T')}.$
%Is $\mathcal T'$ with (RDP)???%

If $s$ is a $\sigma$-additive state on $\mathcal T,$ then
$$s(f)= \int_X f \dx \mu_s, \ f \in \mathcal T',\eqno(4.9)
$$
where $\mu_s(A):= s(\chi_A),$ $A \in \mathcal B_0\mathcal {(T)}.$
\end{proposition}

\begin{proof}  Due to Proposition \ref{pr:8.1}, $\mathcal B_0\mathcal {(T)}$ is
a $\sigma$-algebra of subsets. Hence, $\mathcal T'$ is an
effect-tribe that is a subalgebra of $\mathcal T.$  If $f=\chi_A \in
B(\mathcal T),$ then $\chi_A\in \mathcal T'.$ Now if  $g \in
\mathcal T',$  due to Lemma \ref{le:9.1'}, $gf\in \mathcal T$ so
that $gf \in \mathcal T',$ $g\wedge_{\mathcal T'}
f=gf=g\wedge_{\mathcal T} f.$  Let $g = g_1+g_2,$ where $g_1,g_2\in
\mathcal T$ and $g_1\le \chi_A,$ $g_2\le 1-\chi_A,$ we have
$g_1=g\chi_A \in \mathcal T'$ and $g_2=g(1-\chi_A)\in \mathcal T'.$
This entails $f \in B(\mathcal T')$ and $\mathcal B_0(\mathcal
T)\subseteq \mathcal B_0(\mathcal T').$ On the other hand, if
$f=\chi_A \in B(\mathcal T'),$ then $f\wedge (1-f)=0$ for $\wedge$
taken in $\mathcal T$ as well as in $\mathcal T',$ so that $f \in
B(\mathcal T).$  Hence, $\mathcal B_0\mathcal {(T)}=\mathcal B_0
\mathcal{(T')}.$

Using the proof of Theorem \ref{th:8.3}, we have (4.9).
\end{proof}

We note that in the case of Example \ref{ex:8.2}, $\mathcal T'$
consists from all functions $f \in \mathcal T$ such that $f(a)=f(b)$
and, for each $\sigma$-additive state $s$ on $\mathcal T$ and each
$f \in \mathcal T',$ (4.9) holds.

Now we apply the Butnariu--Klement Theorem for $\sigma$-additive
states on a $\sigma$-complete MV-algebra.

\begin{theorem}\label{th:BK}
Let $s$ be a $\sigma$-additive state on a $\sigma$-complete
MV-algebra $M.$  Then there are a tribe $\mathcal T$ of functions on
some $X\ne \emptyset,$ a $\sigma$-MV-homomorphism $h$ from $\mathcal
T$ onto $M$ and a unique $\sigma$-additive probability measure
$\mu_s$ on $\mathcal B_0\mathcal {(T)}=\mathcal S_0(\mathcal T)$
such that
$$
s(a)=\int_X f(x)\dx \mu_s(x),\quad a \in M, \eqno(4.10)
$$
where $f \in \mathcal T$ and $h(f)=a.$
\end{theorem}

\begin{proof}  By the Loomis--Sikorski Theorem for $\sigma$-complete
MV-algebras, \cite{Dvu3, Mun}, there is a tribe $\mathcal T$ of
functions on some nonempty set $X$ and a $\sigma$-MV-homomorphism
$h$ from $\mathcal T$ onto $M$.  Define $\mathcal B_0\mathcal
{(T)},$ then $\mathcal B_0\mathcal {(T)}=\mathcal S_0(\mathcal T)$
and $\mathcal B_0\mathcal {(T)}$ is a $\sigma$-algebra of subsets of
$X.$ At any rate, $h(\mathcal B_0\mathcal {(T)})\subseteq B(M).$

Define $\mu_s(A)= s(h(\chi_A)),$ $A \in \mathcal B_0\mathcal {(T)},$
and $\hat s(f)=s(h(f)),$ $f\in \mathcal T.$ Then $\mu_s$ is a
$\sigma$-additive probability measure on $\mathcal B_0\mathcal
{(T)}$ and $\hat s$ is a $\sigma$-additive measure on the tribe
$\mathcal T.$  By (4.4), we have
$$ \hat s(f)= \int_X f(x)\dx \mu_s(x),\quad f \in
\mathcal T.\eqno(4.11)
$$
Now if $a =h(f)\in M$ for some $f \in \mathcal T,$ then (4.11)
implies (4.10).

Let $\mu$ be any $\sigma$-additive probability measure on $\mathcal
B_0(\mathcal T)$ such that (4.10) holds for $\mu.$ If $A\in \mathcal
B_0(\mathcal T),$ then $\mu_s(A)= s(h(\chi_A))=\int \chi_A\dx \mu
=\mu(A).$
\end{proof}

\begin{remark}\label{re:BK} {\rm

(1) Theorem \ref{th:BK} holds for any tribe $\mathcal T$ and any
$\sigma$-MV-homo\-morphism $h$ from $\mathcal T$ onto $M.$

 (2) In Theorem \ref{th:BK} we can put $X =\partial_e \mathcal
S(M)$ and by \cite{Dvu3} if  for any $a \in M$, we define a function
$\hat a:\partial_e \mathcal S(M) \to [0,1],$ then $\mathcal T$ is a
tribe generated by $\{\hat a: a\in M\}$ and  $h(f)=a$ iff $\{x\in X:
f(x)=\hat a(x)\}$ is a meager set. Then $h(\mathcal B_0\mathcal
{(T)})=B(M).$

}
\end{remark}

\section{Monotone $\sigma$-Complete Effect Algebras}%5

In the present section, we characterize monotone $\sigma$-complete
effect algebras analyzing their simplex structure of state spaces.
We describe also Jauch-Piron states by the $\sigma$-convex hull of
extremal states.

We start with some necessary definitions on convex structures. For a
more detailed study we recommend \cite{Goo, Alf, AlSc}.

Let $K$ be a compact convex subset of a locally convex Hausdorff
space. A mapping $f:\ K \to \mathbb R$ is said to be {\it affine}
if, for all $x,y \in K$ and any $\lambda \in [0,1]$, we have
$f(\lambda x +(1-\lambda )y) = \lambda f(x) +(1-\lambda ) f(y)$. Let
$\Aff(K)$ be the set of all continuous affine functions on $K.$ Then
$\mbox{Aff}(K)$ is a unital po-group with the strong unit $1$ which
is a subgroup  of the po-group $\mbox{C}(K)$ of all continuous
real-valued functions on $K$ (we recall that, for $f,g \in
\mbox{C}(K),$ $f \le g$ iff $f(x)\le g(x)$ for any $x \in K$), hence
it is an Archimedean unital po-group with the strong unit $1$. In
addition, $\mbox{C}(K)$ is  an $\ell$-group (= lattice ordered
group).

Let $S = {\mathcal S}(\mbox{Aff}(K),1).$ Then the evaluation mapping
$\psi:\ K \to S$ defined by $\psi(x)(f)=f(x)$ for all $f \in
\mbox{Aff}(K)$ $(x \in K)$ is an affine homeomorphism of $K$ onto
$S,$ see \cite[Thm 7.1]{Goo}.

We recall that a {\it convex cone} in a real linear space $V$ is any
subset $C$ of  $V$ such that (i) $0\in C,$ (ii) if $x_1,x_2 \in C,$
then $\alpha_1x_1 +\alpha_2 x_2 \in C$ for any $\alpha_1,\alpha_2
\in \mathbb R^+.$  A {\it strict cone} is any convex cone $C$ such
that $C\cap -C =\{0\},$ where $-C=\{-x:\ x \in C\}.$ A {\it base}
for a convex cone $C$ is any convex subset $K$ of $C$ such that
every non-zero element $y \in C$ may be uniquely expressed in the
form $y = \alpha x$ for some $\alpha \in \mathbb R^+$ and some $x
\in K.$

We recall that in view of \cite[Prop 10.2]{Goo}, if $K$ is a
non-void convex subset of $V$ and if we set
$$
C =\{\alpha x:\ \alpha \in \mathbb R^+,\ x \in K\},
$$
then $C$ is a convex cone in $V,$ and $K$ is a base for $C$ iff
there is a linear functional $f$ on $V$ such that $f(K) = 1$ iff $K$
is contained in a hyperplane in $V$ which misses the origin.

Any strict cone $C$ of $V$ defines a partial order $\le_C$ via $x
\le_C y$ iff $y-x \in C.$ It is clear that $C=\{x \in V:\ 0 \le_C
x\}.$ A {\it lattice cone} is any strict convex cone $C$ in $V$ such
that $C$ is a lattice under $\le_C.$

A {\it simplex} in a linear space $V$ is any convex subset $K$ of
$V$ that is affinely isomorphic to a base for a lattice cone in some
real linear space. A  simplex $K$ in a locally convex Hausdorff
space is said to be (i) {\it Choquet} if $K$ is compact, and (ii)
{\it Bauer} if $K$ and $\partial_e K$ are compact.

Choquet and Bauer simplices can be characterize as follows: (i) if
$E$ is with (RDP), then $\mathcal S(E)$ is a Choquet simplex,
\cite[Thm 10.17]{Goo}. Let $K$ be a convex compact subset of a
locally convex Hausdorff space, then (ii) $K$ is a Choquet simplex
iff $(\mbox{Aff}(K),1)$ is an interpolation po-group, \cite[Thm
11.4]{Goo}, (iii) $\mathcal S(E)$ is a Bauer simplex whenever $E$ is
an MV-algebra, (iv) $K$ is a Bauer simplex iff $(\mbox{Aff}(K),1)$
is an $\ell$-group, \cite[Thm 11.21]{Goo}. (v) The state space of
$\mathcal E(H)$ is not a simplex, \cite[Ex 4.2.6]{BrRo}.

We say that an effect algebra $E$ satisfies {\it general
comparability} if, given $x,y \in E$, there is a central element $e
\in E$ such that $p_e(x) \le p_e(y)$ and $p_{e'}(x) \ge p_{e'}(y)$
where $p_e(x) = x\wedge e.$ This means that the coordinates of the
elements $x = (p_e(x),p_{e'}(x))$ and $y=(p_e(y),p_{e'}(y))$ can be
compared in $[0,e]$ and $[0,e']$, respectively. If $E$ satisfies
(RDP) and general comparability, then $E$ is a lattice, and it can
be converted  to an MV-algebra where original $+$ and that derived
from the MV-structure coincide.

For example, (i) every linearly ordered pseudo-effect algebra
trivially satisfies general comparability; (ii) also any Cartesian
product of linearly ordered pseudo-effect algebras. If $E$ satisfies
general comparability, then $E$ is a lattice with (RDP) and it can
be organized into an MV-algebra.  In addition, every
$\sigma$-complete MV-algebra satisfies general comparability,
\cite[Thm 9.9]{Goo}. Moreover, if $E$ satisfies general
comparability, every extremal state on $B(E)$ can be uniquely
extended to an extremal state on $E,$ \cite[Thm 8.14]{Goo},
\cite[Thm 4.6]{Dvu2}.

In Example \ref{ex:8.2}, we have a case that $\mathcal T$ is an
effect-tribe with (RDP)  such that $\mu_0$ is an extremal state on
$B(E)$ but it has an extension to a state on $\mathcal T$ that is
not extremal; moreover both $\mu_0$ and $\mu_1$ have the same
restriction to $B(\mathcal T)$.

In what follows, we show when the extension is unique.

We assert that if $s$ is an extremal state on $E,$ then its
restriction, $s_B= s\,\upharpoonright\, B(E),$ to $B(E)$ is also an
extremal state on $B(E).$ Indeed, take $e\in B(E)$, we assert that
$s(e)\in \{0,1\}.$ If not, then $0<s(e)<1$ and  $s= \lambda s_1
+(1-\lambda)s_2$ where $\lambda = s(e),$ $s_1(x)=s(x\wedge e)/s(e)$
and $s_2(x)= s(x\wedge e')/s(e')$ which gives a contradiction.
Therefore, $s_B$ is an extremal state on $B(E),$  and the mapping
$\theta:\ \partial_e {\mathcal S}(E)\to
\partial_e {\mathcal S}(B(E))$ given by
$$ \theta(s) = s_B := s \upharpoonright B(E),\ s \in
\partial_e{\mathcal S}(E), \eqno(5.1)
$$
is well defined and continuous.

\begin{theorem}\label{th:8.6} Let $E$ be a monotone
$\sigma$-complete effect algebras with {\rm (RDP)}. The following
statements are equivalent:

\begin{enumerate}
\item[{\rm (i)}] Every extremal state on $B(E)$ is extendible to a unique
state on $E$ and this state is extremal.

\item[{\rm (ii)}]  ${\mathcal S}(E)$ is a Bauer simplex.

\item[{\rm (iii)}] $E$ is lattice ordered.

\item[{\rm (iv)}]  $E$ satisfies general comparability.

\item[{\rm (v)}] The mapping  $\theta:\ \partial_e {\mathcal S}(E)\to
\partial_e {\mathcal S}(B(E))$ defined by {\rm (5.1)} is a
homeomorphism.

\item[{\rm (vi)}]  $E$ can be converted into an MV-algebra, where
original $+$ and derived one from $\oplus$ coincide.
\end{enumerate}

\end{theorem}

\begin{proof}  (i) $\Rightarrow$ (ii). Let $\{s_\alpha\}$ be a net of
extremal states on $E$ converging to a state $s \in {\mathcal
S}(E)$. We assert that $s$ is extremal. It is clear that $s_\alpha
\upharpoonright B(E)$ converges weakly to $s_B :=s\upharpoonright
B(E)$. But the space $\partial_e{\mathcal S}(B(E))$ is compact, so
that $s_B \in \partial_e{\mathcal S}(B(E)).$ Because $s_B$ has a
unique extension, $s$, therefore, $s$ is extremal, and ${\mathcal
S}(E)$ is a Bauer simplex.

According to \cite[Cor 16.28]{Goo}, (ii)--(iv) are equivalent.

(iv) $\Rightarrow$ (i). This follows from \cite[Prop 8.13]{Goo} or
\cite[Thm 4.4, Cor 4.5]{Dvu2}.

(iv) $\Rightarrow$ (v). This follows from \cite[Thm 8.14]{Goo} and

(v) $\Rightarrow$ (ii).  This is evident.

(iv) $\Rightarrow$ (vi). This was already mentioned.

(vi) $\Rightarrow$ (iii). It is evident.
 \end{proof}

Now we reduce  criterion (3.1) for extremality of states for
monotone $\sigma$-complete effect algebras satisfying (RDP).

\begin{proposition}\label{pr:8.7}  Let $E=\Gamma(G,u)$ be a monotone
$\sigma$-complete effect algebra satisfying {\rm (RDP)}. A state $s$
on $(G,u)$ is extremal if and only if, given $f,g \in G^+,$ there is
$h \in G^+$ such that $h\le f,g$ and $\min\{s(f),s(g)\} = s(h).$
\end{proposition}

\begin{proof} If the criterion is satisfied, then $s$ is extremal.
Converse, let $s$ be an extremal state on $(G,u).$ Denote $s_0 =
\sup\{s(h): h \in G^+,\ \ h\le f,g\}=\min\{s(h),s(g)\}.$ There
exists a sequence of elements $\{h_n\}$ in $(G,u)$ such that $h_n\le
f,g$ and $\lim_n s(h_n) = s_0.$ Due the the interpolation property
holding in the po-group $G,$ there exists a monotone sequence
$\{h_n'\}$ such that, for any $n\ge 1,$ $h_1,\ldots,h_n \le h_n' \le
f,g.$  Let $h_0 = \bigvee _n h_n' \in G,$ then $h_n \le h_0 \le f,g$
so that $s(h_n) \le s(h_0) \le s_0$ and $s(h_0)=s_0.$
\end{proof}

Proposition \ref{pr:8.7} can be generalized as follows.

We recall that a poset $E$ satisfies the {\it countable
interpolation property} provided that for any two sequences of
elements of $E$, $\{a_i\}$ and $\{b_j\},$ such that $a_i\le b_j$ for
all $i,j,$ there is an element $c\in E$ such that $a_i\le c\le b_j$
for all $i,j.$ Due to \cite[Prop 16.3]{Goo}, an effect algebra
$E=\Gamma(G,u)$ with (RDP) has countable interpolation  iff the
po-group $G$ has countable interpolation, and due to \cite[Thm
16.10]{Goo}, if an effect algebra $E$ with (RDP) is monotone
$\sigma$-complete, then $E$ satisfies countable interpolation.

On the other hand, there is even an MV-algebra satisfying countable
interpolation that is not monotone $\sigma$-complete as we can
deduce from \cite[p. 280]{Goo}.

\begin{proposition}\label{pr:8.7.1}
Let $E=\Gamma(G,u)$ be an  effect algebra satisfying {\rm (RDP)} and
countable interpolation. A state $s$ on $(G,u)$ is extremal if and
only if, given $f,g \in G^+,$ there is $h \in G^+$ such that $h\le
f,g$ and
$$s(h)=\min\{s(f),s(g)\}.\eqno(5.2)$$

If, in addition, $E$ is divisible, then a state $s$ on $E$ is
extremal if and only if {\rm (5.2)} holds for all $f,g \in E.$
\end{proposition}

\begin{proof} One implication is evident.

Now let $s$ be an extremal state on $(G,u).$ Denote $s_0 =
\sup\{s(h): h \in G^+,\ \ h\le f,g\}=\min\{s(h),s(g)\}.$ There
exists a sequence of elements $\{h_n\}$ in $(G,u)$ such that $h_n\le
f,g$ and $\lim_n s(h_n) = s_0.$ Due the countable interpolation
property holding in the po-group $G,$ there exists an element
$h_0\in G$ such that $h_n\le h_0 \le f,g,$ for each $n,$ so that
$s(h_n) \le s(h_0) \le s_0$ and $s(h_0)=s_0.$

The last statement follows from Proposition \ref{pr:3.2}.
\end{proof}

The latter two propositions can be used to characterize Jauch-Piron
states. We say that a state $s$ on an effect algebra $E$ is {\it
Jauch-Piron} if, for all $a,b\in E$ with $s(a)=s(b)=1,$ there is an
element $c\in E$ such that $c\le a,$ $c\le b$ and $s(c)=1.$  For
example, if $E$ is an MV-algebra, then every state $s$ is
Jauch-Piron in view of the property $s(a\vee b)+s(a\wedge
b)=s(a)+s(b)$ holding in each MV-algebra. Every state in Example
\ref{ex:4.1} is also Jauch-Piron.

We recall that if $\{s_n\}$ is a finite or a countable set of
states, then $s=\sum_n\lambda_ns_n,$ where $\lambda_n \ge 0$ and
$\sum_n\lambda_n=1,$ is a state that is called a $\sigma$-{\it
convex combination} of $\{s_n\}.$

\begin{proposition}\label{pr:8.7.2}
Let $E$ be an  effect algebra satisfying {\rm (RDP)} and countable
interpolation. Then every $\sigma$-convex combination of extremal
states is Jauch-Piron.
\end{proposition}

\begin{proof}  Suppose that $s$ is an extremal state on $E$ and let
$s(a)=s(b)=1.$  Due to Proposition \ref{pr:8.7.1}, there is an
element $c\in E$ such that $c\le a,b$ and
$s(c)=\min\{s(a),s(b)\}=1.$

Now let $s=\sum_n \lambda_ns_n$ where each $s_n\in \partial_e
\mathcal S(E)$ and $0<\lambda_n<1,$ $\sum_n \lambda_n=1.$ Let
$s(a)=s(b)=1.$ Then $s_n(a)=s_n(b)=1$ for each $n.$  Hence, there is
$c_n\le a,b$ such that $s_n(c_n)=1.$ Due to the countable
interpolation property, there is an element $c\in E$ such that $c_n
\le c\le a,b$ for every $n.$ Then $s_n(c_n)\le s_n(c)=1$ and
$s(c)=1.$
\end{proof}

Now we present some examples to show how the extremality criterion
(5.2) works  on effect-clans ${\mathcal T} \subseteq [0,1]^X$ with
(RDP) for states of the form
$$s_x(f)=f(x),\ f\in \mathcal T, \eqno(5.3)
$$
where $x \in X.$ We note that each $s_x$ is in fact a
$\sigma$-additive state on the effect-tribe $\mathcal T$ and the
system $\{s_x: x \in X\}$ is an order determining system of states.

\begin{example}\label{ex:8.8}  {\rm (1)  Let $X=[0,1]$ be the real
unit interval and $E=\Gamma(\mbox{C}(X),1_X)).$ Since $X$ is a
Hausdorff compact topological space, due to the Riesz Representation
Theorem, a state $s$ is extremal iff $s=s_x$ for some $x \in X,$
moreover, if $s_x(f)=\int_X f \dx \mu$ where $\mu$ is the Borel
probability on the Borel $\sigma$-algebra ${\mathcal B}_0(X),$ then
$\mu = \delta_x.$  The same is true if $X$ is any compact Hausdorff
topological space and $\mu$ is a regular Borel probability measure
${\mathcal B}_0(X).$

\vspace{3mm}\noindent (2) Let $X=[0,1]$ be the real interval. Then
$X$ is a convex compact Hausdorff topological space; set $E =
\Gamma(\mbox{Aff}(X),1_X).$  Because $X$ is a simplex (for the
definition of a simplex see \cite{Goo} or the next chapter), $E$ is
an effect-clan with (RDP), \cite[Thm 11.4]{Goo}. Every element $f
\in \mbox{Aff}(X)$ is of the form $f(x)=\alpha x +\beta$ where $0\le
\beta\le 1$ and either $0\le \alpha\le 1$ or $-1\le \alpha <0.$ $E$
is an effect-tribe: Indeed, let $f_n(x) =\alpha_nx+\beta_n\le
f_{n+1}(x)=\alpha_{n+1}x+\beta_{n+1}$ be bounded in $E$. Then
$\beta_n\le \beta_{n+1}$ so that there exists $\beta_0 = \lim_n
\beta_n$ and due to $\alpha_n = (f_n(x)-\beta_n)/x$ for $x
>0,$ there exists $\lim_n \alpha_n = \alpha_0.$  Let $n,k$ be
arbitrary integers, then $f_n(x)\le
f_{n+k}(x)=\alpha_{n+k}x+\beta_{n+k},$ so that $f_n(x)\le
\lim_kf_{n+k}(x)=\alpha_0x+\beta_0.$ Hence, $f(x)=\alpha_0
x+\beta_0\in E.$  Moreover, $E$ is a lattice such that $f\wedge g
\in E$ (with respect to the order by points) for every $f,g \in E,$
but $\min\{f,g\}$ is not necessarily in $E$, so that $E$ is not a
tribe.

We state that  $s_0$ and $s_1$ are extremal states, but $s_x
\notin\partial_e{\mathcal S}(E)$ for any $x \in (0,1):$   Since $X$
and $\mathcal S(E)$ are affinely homeomorphic,  $\partial_e
{\mathcal S}(E)=\{s_0,s_1\}$ and ${\mathcal S}(E)=\{s_x:\ x \in
[0,1]\}.$

\vspace{3mm}\noindent (3)   Let the effect-tribe be that from
Example \ref{ex:8.2}. Then every state $s_x$ for $x \in X$ is
extremal, where $s_x(f)=f(x),$ $f \in {\mathcal T}.$ Indeed, let
$x_0$ and $f,g \in G^+$ be fixed. We verify criterion (5.2).

(i)  Assume that $x_0 \in X\setminus \{a,b\}.$ Since $f,g$ are
bounded, take $h\in G^+$ such that $h(x)=0$ for each $x\ne x_0$ and
$h(x_0) = \min\{f(x_0),g(x_0)\}.$ By (5.2), $s_{x_0}$ is extremal.

(ii)  Assume  that $x_0 =a$ and choose  $h\in G^+$ such that $h(a) =
\min\{f(a),g(a)\},$ $h(b) =0,$ $h(x) =(h(a)+h(b))/2=h(a)/2$ for $x
\in (A_f\cap A_g)\setminus \{a,b\}$ and $h(x)=0$ for otherwise.
Again by (5.2), $s_a$ is extremal and the same true for $s_b.$

Nevertheless $\mu_0$ is extremal, the state $s_{\mu_0}$ is not
extremal because $s_{\mu_0}= (s_a+s_b)/2$ but in view of Proposition
\ref{pr:8.7.2}, $s_{\mu_0}$ is a Jauch-Piron state.

\vspace{3mm}\noindent (4) If we consider the tribe $\mathcal T(H),$
then every state $s_x(f_A)=(Ax,x)$ is an extremal state when $\dim
H\ge 3.$ This follows from the Aarnes theorem and not from (5.2).
because the state space of $\mathcal E(H)$ is not a simplex.

}
\end{example}

\section{States on Effect Algebras and Integrals}

In the present section, we show that any state on an effect algebra
can be represented as an integral through a regular Borel
probability measure. This will generalize the results from
\cite{Dvu4} where this was proved only for states on interval effect
algebras.

We start with some definitions.

If  $K$ is a compact Hausdorff topological space, let ${\mathcal
B}(K)$ be the Borel $\sigma$-algebra of $K$ generated by all open
subsets of $K.$  Let  ${\mathcal M}_1^+(K)$ denote the set of  all
probability measures, that is, all positive regular
$\sigma$-additive Borel measures $\mu$ on $\mathcal B(K).$  We
recall that a Borel measure $\mu$ is called regular if

$$\inf\{\mu(O):\ Y \subseteq O,\ O\ \mbox{open}\}=\mu(Y)
=\sup\{\mu(C):\ C \subseteq Y,\ C\ \mbox{closed}\}
$$
for any $Y \in {\mathcal B}(K).$

The following result is motivated by research in \cite{KuMu}. Here
we present another proof which is motivated by the theory of states
on effect algebras.

Let $E\ne \emptyset$ and let $[0,1]^E$ be endowed with a product
topology, i.e., with the weak topology, and let $\mathcal S$ be a
nonempty closed convex subset of the Tikhonov cub $[0,1]^E.$ Any map
$s \in \mathcal S$  will called also a state. For any $a \in E,$ let
$\hat a(s):= s(a),$ $s \in S.$ Then $\widehat E=\{\hat a: a\in E\}
\subseteq \Aff(\mathcal S).$

\begin{theorem}\label{th:ps2}  Let $s \in \mathcal S$ and let $\mathcal S$ be
convex and closed in $[0,1]^E.$  Then there is a regular Borel
probability measure, $\mu_s,$ on $\mathcal {B(S)}$ such that

$$s(a)=\int_{\mathcal S} \hat a(x)\dx \mu_s(x),  \eqno(6.1)
$$
for each $a \in E.$
\end{theorem}

\begin{proof} We define a congruence $\equiv$ on $E$ by $a \equiv b$
iff $s(a)=s(b)$ for any $a \in E.$  Given $s \in \mathcal S,$ denote
by $\hat s$ the mapping from $E/\equiv$ into  $[0,1]$ defined by
$\hat s([a])=s(a),$ where $[a]$ denotes the congruence class
corresponding to the element $a \in E,$ and let $\widehat {\mathcal
S}=\{\hat s: s \in \mathcal S\}.$ Then the mapping $s \mapsto \hat
s$ is an affine homeomorphism.  Therefore, without loss of
generality, we can assume that if $a\ne b,$ then there is a state $s
\in \mathcal S$ such that $s(a)\ne s(b).$

Denote by $G$ the subgroup of $\Aff(\mathcal S)$ generated by
$\widehat E,$ and let $E_1 =\Gamma(G,1)$ be the effect algebra
generated by $\widehat E.$

We have that $\mathcal S$ and $\mathcal S(\Aff(\mathcal S),1)$ are
affinely homeomorphic under the mapping $\rho(s)(f)= f(s),$ $ f \in
\Aff(\mathcal S).$  Hence, the restriction of $\rho (s)$ onto $E_1$
is a state of the effect algebra $E_1.$  Now let $t$ be a state on
$E_1.$ By \cite[Cor 4.3]{Goo} or by \cite[Thm 6.9]{Dvu4}, this state
can be extended to a state $\hat t$ on $\Gamma(\Aff(\mathcal S),1)$
and to a state on $(\Aff(\mathcal S),1).$  Hence, there is a mapping
$s\in \mathcal S$ such that $\hat t = \rho(s),$ so that $t(g)=g(s)$
for any $g \in G.$ Due to our assumption, we see that if
$\rho(s_1)=\hat t=\rho(s_2),$ then $s_1 = s_2.$ Therefore, the
mapping $t\in \mathcal S(E_1)\mapsto s,$ where $s \in \mathcal S$
such that $\rho(s)|E_1=t$ is an affine homeomorphism; and  we denote
by $\theta:\mathcal S\to \mathcal S(E_1) $ the inverse affine
homeomorphism.

So that, if $s \in \mathcal S$ and $a\in E,$ then $s(a) = \hat a(s)
=\rho(s)(\hat a)$ and using \cite[Thm 6.3]{Dvu4}, there is a regular
Borel probability measure $\nu_s$ defined on $\mathcal B(\mathcal
S(E_1))$ such that
$$s(a)=\rho(s)(\hat a) = \int_{\mathcal S(E_1)} \widehat{\hat
a}(y)\dx\nu_s(y),
$$
where $\widehat{\hat a}$ is a mapping from $\mathcal S(E_1)$ into
$[0,1]$ such that $\widehat{\hat a}(t):=t(\hat a),$ $t \in \mathcal
S(E_1).$  If we define a regular Borel measure $\mu_s=\nu_s\circ
\theta$ defined on $\mathcal {B(S)},$ we have

$$s(a)= \int_{\mathcal \theta^{-1}(\mathcal S(E_1))} \widehat{\hat
a}(\theta(x))\dx\nu_s(\theta(x)) = \int_{\mathcal S} \hat
a(x)\dx\mu_s(x).
$$
\end{proof}

\begin{remark}\label{re:ps1} {\rm If $\mathcal S$ is a Choquet simplex,
then a regular Borel probability measure on $\mathcal {B(S)}$  in
(6.1) is a unique maximal regular Borel probability measure such
that $\mu_s \sim \delta_s$ (see \cite[Thm 6.3]{Dvu4}); if $\mathcal
S$ is a Bauer simplex, then $\mu_s$ is a unique regular Borel
probability measure, and $\mu_s(\partial_e\mathcal S)=1.$}
\end{remark}

In what follows, we show that (6.1) holds also for any nonempty
compact set $\mathcal S \subseteq [0,1].$ We note that we will prove
it using other argumentations than those used in Theorem
\ref{th:ps2}.

\begin{theorem}\label{th:ps3}  Let $s \in \mathcal S$ and let $\mathcal S$ be
closed in $[0,1]^E.$  Then there is a regular Borel probability
measure, $\mu_s,$ on $\mathcal {B(S)}$ such that {\rm (6.1)} holds
for any $a \in E$ and $s \in \mathcal S.$
\end{theorem}

\begin{proof}  The mapping $\epsilon: \mathcal S \to K:=\mathcal
{M}^+_1(\mathcal S)$ defined by  $\epsilon(s) =\delta_s,$ where
$\delta_s$ is the Dirac measure concentrated at the point $s$, is a
homeomorphism from $\mathcal S$ onto $\partial_e \mathcal
{M}^+_1(\mathcal S).$

Let $a \in E$, and let $a_\epsilon: \partial_e K \to[0,1]$ be
defined by $a_\epsilon(\delta_s):=s(a).$  Then $a_\epsilon \in
\mbox{C}(\partial_e K).$ Since $K$ is a simplex, \cite[Cor
10.18]{Goo}, by the Tietze Theorem \cite[Prop II.3.13]{Alf},
$a_\epsilon$ can be uniquely extended to an affine function $\hat
a_\epsilon\in \Gamma(\Aff(K),1)$ defined on $K.$  If $\hat
a:\mathcal S\to [0,1],$ where $\hat a(s):=s(a),$ then $\hat
a_e(\epsilon (s))=\hat a(s)$ for each $s\in \mathcal S.$

Since the effect algebra $\Gamma(\mbox{C}(K),1)$ is an MV-algebra,
\cite[Thm 6.5]{Dvu4} yields there exists a (unique) regular Borel
probability measure $\nu_s$ on $\mathcal B(K)$ such that
$\nu_s(\partial_e K)=1$ and
$$
s(a)=\hat a_\epsilon(\delta_s)=\int_{\partial_e K} \hat
a_\epsilon(y)\dx\nu_s(y).
$$
If we set $\mu_s:= \nu_s\circ \epsilon$, then $\mu_s$ is a regular
Borel probability measure on $\mathcal {B(S)}.$ Therefore,
$$
s(a)= \int_{\epsilon^{-1}(\partial_e K)} \hat a_\epsilon
(\epsilon(x))\dx \nu_s(\epsilon(x))= \int_{\mathcal S} \hat a(x)\dx
\mu_s(x)
$$
that finishes the proof.
\end{proof}

\begin{corollary}\label{co:int}  Let $s$ be a state on an effect
algebra $E.$ There is a regular Borel probability measure, $\mu_s$,
on the Borel $\sigma$-algebra $\mathcal B(\mathcal S(E))$ such that
$$
s(a) = \int_{\mathcal S(E)} \hat a(x)\dx\mu_s(x), \quad a\in E.$$

If, in addition, $\mathcal S(E)$ is a Choquet simplex, then there is
a unique regular Borel probability measure, $\mu_s,$ on $\mathcal
B(\mathcal S(E))$ such that $\mu_s(\partial_e \mathcal S(E))=1$ and

$$
s(a) = \int_{\partial_e \mathcal S(E)} \hat a(x)\dx\mu_s(x), \quad
a\in E.$$

\end{corollary}

\begin{proof}  It follows from Theorem \ref{th:ps2} and Remark
\ref{re:ps1}.
\end{proof}

\end{document}